\documentclass[10pt]{amsart}
\usepackage{amsmath}
\usepackage{latexsym}
\usepackage{amssymb,amsthm,amsfonts}
\usepackage{graphicx}
\usepackage[active]{srcltx}
\usepackage{color} 
\usepackage{esint}

\renewcommand{\a }{\alpha }
\renewcommand{\b }{\beta }
\renewcommand{\d }{\delta }
\newcommand{\D }{\Delta }

\newcommand{\e }{\varepsilon }
\newcommand{\g }{\gamma}

\renewcommand{\l }{\lambda }
\renewcommand{\L }{\Lambda }

\newcommand{\n }{\nabla }
\newcommand{\var }{\varphi }

\newcommand{\Sig }{\Sigma}

\renewcommand{\th }{\theta }

\newcommand{\ov}{\overline}
\newcommand{\intbar}{\mathop{\int\makebox(-13.5,0){\rule[4pt]{.7em}{0.3pt}}%
\kern-6pt}\nolimits}

\newcommand{\be}{\begin{equation}}
\newcommand{\ee}{\end{equation}}
\newenvironment{pf}{\noindent{\sc Proof}.\enspace}{\rule{2mm}{2mm}\medskip}

\newcommand{\R}{\mathbb{R}}

\newcommand{\N}{\mathbb{N}}
\newcommand{\pa}{\partial}
\newcommand{\dis}{\displaystyle}
\newcommand{\cD}{\mathcal{D}}
\newcommand{\cJ}{\mathcal{J}}

\newtheorem{lem}{Lemma}[section]
\newtheorem{pro}[lem]{Proposition}
\newtheorem{thm}[lem]{Theorem}
\newtheorem{rem}[lem]{Remark}
\newtheorem{cor}[lem]{Corollary}

\begin{document}

\title[New improved Moser-Trudinger inequalities] {New improved Moser-Trudinger
inequalities and Singular Liouville equations on compact surfaces}

\author{Andrea Malchiodi and David Ruiz}

\address{SISSA, via Bonomea 265, 34136 Trieste (Italy) and Departamento de
An\'alisis Matem\'atico, University of Granada, 18071 Granada
(Spain).}

\email{malchiod@sissa.it, daruiz@ugr.es}

\keywords{Geometric PDEs, Variational Methods, Min-max Schemes.}

\subjclass[2000]{35B33, 35J35, 53A30, 53C21.}

\begin{abstract}
We consider a singular Liouville equation on a compact surface,
arising from the study of Chern-Simons vortices in a self dual
regime. Using new improved versions of the Moser-Trudinger
inequalities (whose main feature is to be scaling invariant) and a
variational scheme, we prove new existence results.
\end{abstract}

\maketitle

\section{Introduction}\label{s:in}

\noindent We consider a compact orientable surface $\Sig$ with
metric $g$ and the equation
\begin{equation}\label{eq:e-1}
    - \Delta_g u = \rho \, h(x)e^{2u}
    - 2 \pi \sum_{j=1}^m \alpha_j  \delta_{p_j} + c, \quad \int_{\Sig}
    h(x)e^{2u} d V_g=1.
\end{equation}

Here $\rho$ is a positive parameter, $h : \Sig \to \R$ a smooth
positive function, $\alpha_j \in [0,1]$, $p_j \in \Sig$ and $c$ is
a constant. Integrating by parts we get that a necessary condition
for the existence of solution is $c= \left (2 \pi \sum_{j=1}^m
\alpha_j - \rho \right)|\Sig|^{-1}$.

This equation arises from physical models such as the abelian
Chern-Simons-Higgs theory and the Electroweak theory, see
\cite{dunne}, \cite{hkp}, \cite{jawe}, \cite{lai}. We also refer
to \cite{tar}, \cite{tardcds}, \cite{yang} and the bibliographies
therein for a more recent and complete description of the subject.
Here we limit ourselves to mention that $u$ is related to the
absolute value of the wave function in the above models, while the
$p_j$'s, called {\em vortices}, are points where the wave function
vanishes. Equation \eqref{eq:e-1} has been the subject of several
investigations, see for example \cite{barsi1}, \cite{barsi2},
\cite{bclt}, \cite{barlin}, \cite{barlintar}, \cite{btjde02},
\cite{btcmp02}, \cite{chakie}, \cite{chouwan}, \cite{det},
\cite{linwang},  \cite{pratar}, \cite{tjfa05}, \cite{tind05},
\cite{luzh2}.

Most of these results deal with asymptotic analysis or compactness
of solutions, while relatively few ones are available concerning
existence. In \cite{dpem}, \cite{esp} some perturbative results
are given, providing solutions of multi-bump type for special
values of the parameter $\rho$. In \cite{btcmp02} an existence
theorem is proved for surfaces with positive genus and for $\rho
\in (4\pi, 8\pi)$. Finally, in \cite{cl4} the Leray-Schauder
degree is computed for $\alpha \geq 1$ and $\rho \in (4\pi,
8\pi)$, and so existence results are deduced.

Our goal is to develop a global variational theory for the
equation, yielding existence of solutions under rather general
conditions. In this paper we give a new improved Moser-Trudinger
inequality, a basic tool for this strategy, and derive some first
results in this spirit. Further existence results will be
discussed in a forthcoming paper.

\

\noindent It is easy to see that equation \eqref{eq:e-1} is
equivalent to:

\begin{equation}\label{eq:e0}
    - \Delta_g u = \rho \left( \frac{h(x)e^{2u}}{\int_\Sig h(x) e^{2u} d V_g}
    - \frac{1}{|\Sig|} \right) - 2 \pi \sum_{j=1}^m \alpha_j \left( \delta_{p_j} - \frac{1}{|\Sig|}
    \right).
\end{equation}

This last formulation has the advantage that it is invariant with
respect to the addition of constants.

\

\noindent Let $G_p(x)$ be the Green's function of $- \D_g$ on
$\Sig$ with singularity at $p$, namely the unique solution of
$$
   - \D_g G_p(x) = \d_p - \frac{1}{|\Sig|} \quad \hbox{ on }
   \Sig, \ \mbox{ with } \int_{\Sigma} G_p(x) \, dV_g =0.
$$

The change of variables
\begin{equation}\label{eq:change}
    u \mapsto u + 2\pi \sum_{j=1}^m \alpha_j G_{p_j}(x)
\end{equation}
transforms \eqref{eq:e0} into an equation of the form
\begin{equation}\label{eq:ee}
    - \Delta_g u  = \rho \left (\frac{\tilde{h}(x) e^{2u}}{\int_\Sig
    \tilde{h}(x) e^{2u} dV_g} - \frac{1}{|\Sig|} \right ) \qquad \quad \hbox{ on } \Sig.
\end{equation}
Since $G_p$ has the asymptotic behavior $G_p(x) \simeq
\frac{1}{2\pi} \log \frac{1}{d(x,p)}$ near $p$, we have
\begin{equation}\label{eq:tildeh}
    \tilde{h} > 0 \hbox{ on } \Sig \setminus \cup_j \{p_j\}; \qquad \quad
  \tilde{h}(x) \simeq d(x,p_j)^{2 \alpha_j} \quad \hbox{ near }
  p_j.
\end{equation}

\

\noindent Problem \eqref{eq:ee} is the Euler-Lagrange equation for the
functional
\begin{equation}\label{eq:I_rho}
    I_\rho(u) = \int_{\Sig} |\n_g u|^2 dV_g + 2 \frac{\rho}{|\Sig|} \int_\Sig  u
    dV_g - \rho \log \int_\Sig \tilde{h}(x) e^{2u} dV_g; \qquad u \in H^{1}(\Sig).
\end{equation}
Recall the Moser-Trudinger inequality
\begin{equation}\label{eq:mtin}
    \log \int_\Sig e^{2(u - \ov{u})} dV_g \leq \frac{1}{4 \pi}
    \int_\Sig |\n_g u|^2 dV_g + C; \quad \qquad  u \in H^1(\Sig),
\end{equation}
see e.g. \cite{moser}. From that inequality one can easily check
that $I_\rho$ is bounded from below if $\rho < 4 \pi$. Moreover,
$4 \pi$ is a threshold value, in the sense that for larger values
of $\rho$ the functional does not have a finite lower bound.
However one can still hope to find critical points of saddle type,
using for example min-max schemes.

\

\noindent This strategy or other topological methods (jointly with
blow-up estimates), have been used successfully for {\em regular
Liouville equations} of the form \eqref{eq:e0} but with all the
$\a_i$'s equal to zero. Such problems have motivations arising
from physics (study of mean field vorticity, or Chern-Simons
theory without sources), or from conformal geometry (prescribing
the Gauss curvature or some of its higher order counterparts), see
\cite{bls}, \cite{bm}, \cite{calu}, \cite{cygc}, \cite{cl1},
\cite{cl2}, \cite{dem}, \cite{dem2}, \cite{djlw}, \cite{dj},
\cite{dm}, \cite{dr}, \cite{li}, \cite{lisha}, \cite{ll0},
\cite{lucia}, \cite{ns}, \cite{st}. In the latter case, the
function $e^{2u}$ represents the conformal dilation of the
background metric on a given surface. In fact, \eqref{eq:e0} also
arises in the Gauss curvature prescription problem on surfaces
with conical singularities: for more details we refer to
\cite{bdm}.

In this framework, one main common tool for applying variational
arguments is some kind of improvement of the Moser-Trudinger
inequality. A classical example is a result by J.Moser, \cite{mo},
where he showed that the first constant in \eqref{eq:mtin} can be
taken to be $\frac{1}{8 \pi}$ for even functions on $S^2$. A more
general improvement was obtained by T.Aubin in \cite{aubin}, still
in the case of the standard sphere: he showed that for {\em
balanced} metrics one can take any constant which is larger than
$\frac{1}{8 \pi}$ (provided $C$ is taken large enough).
Interesting applications were found for example in \cite{cygc}
where rather general conditions were given for prescribing the
Gauss curvature on the sphere. Aubin's improvement was generalized
by W.Chen and C.Li in \cite{cl} for all surfaces, under the
condition that the conformal volume $e^{2u}$ {\em spreads} into
two distinct regions (separated by a positive distance). This
result was used in \cite{djlw} to produce solutions of the regular
Liouville equation, see also \cite{dj}, \cite{dm} and \cite{mal}
for further progress on this direction.

\

\noindent The main goal of this paper is to obtain a new type of
improved inequality and apply it to the study of \eqref{eq:e0}. To
explain the spirit of this improvement we recall the result in
\cite{tro}, which states that if $\tilde{h}$ is as in
\eqref{eq:tildeh} then the best constant $\mathcal{A}$ for the
inequality
$$
   \log \int_\Sig \tilde{h} e^{2(u-\ov{u})} dV_g \leq \mathcal{A}
   \int_\Sig |\n u|^2 dV_g + C \qquad \left( \ov{u} = \intbar_\Sig u dV_g
   \right)
$$
is given by $\left(4 \pi \min \left\{ 1, \min_i \{(1 +
\a_i)^{-1}\} \right\}\right )^{-1}$. Therefore if some of the
$\a_i$'s is negative the best possible constant is lower than
$\frac{1}{4 \pi}$, but if all the $\a_i$'s are positive (as in our
case) the best constant is just $\frac{1}{4 \pi}$. One can easily
see this by testing the inequality on a {\em standard bubble},
namely a function of the form
\begin{equation}\label{eq:bubble}
    \var_{\l,x}(y) = \log \frac{\l}{1 + \l^2 dist(x,y)^2},
\end{equation}
with center point $x$ different from all the $p_i$'s. This
function realizes the best constant in the regular case, and for
the above choice of $x$ there is basically no effect from the
vanishing of $\tilde{h}$ somewhere on $\Sig$.

On the other hand, in \cite{det} it was shown that for any $\a >
-1$ there exists $C_\a$ such that
$$
  \log \int_B |x|^{2 \a} e^{2(u - \ov{u})} dV_g \leq \frac{1}{4 (1+\a) \pi}
    \int_B |\n_g u|^2 dV_g + C_\a; \qquad \quad u \in H^1_r(B).
$$
In the latter formula $B$ stands for the unit ball of $\R^2$ and
$H^1_r$ denotes the space of radial functions of class $H^1$ in
$B$. Our improvement substitutes the symmetry requirement with a
condition which, heuristically, applies to a subset of functions
having codimension two in $H^1$, assuming $\a \in (0,1]$.
Roughly speaking, we associate to each function $u$ a {\em center
of mass} constructed out of the unit measure $\mu_u :=
\frac{|x|^{2\a} e^{2u}}{\int_B |x|^{2\a} e^{2u} dx}$: then the
improvement would occur for functions whose center of mass is the
origin (in our case, a singular point).

However, a rigorous proof of the above claim requires new
arguments: the proof of Aubin (and in fact also Chen and Li's one)
relies on the fact of being able to find two sets with positive
distance (bounded away from zero) which both contain a finite portion of the total
conformal volume. The positivity of the distance allows to find
cutoff functions $\chi_i$ with bounded gradients and to apply the
standard Moser-Trudinger inequality to $\chi_i u$, choosing then
the $\chi_i u$ with the smaller Dirichlet energy. This strategy
fails in our case since the measure $\mu_u$ can be arbitrarily
concentrated near a single point.

What is needed in this context is a condition which stays
invariant under dilation of the measure $\mu_u$. To achieve it we
use concentration functions (in the spirit of
concentration-compactness results) and covering arguments through
thick annuli, see Section \ref{s:cov} for details. Somehow, we
want to define a {\em continuous} map from $H^1(\Sig)$ into $\Sig$
which keeps track of the points with maximal concentration of
conformal volume: an improvement will hold when the image of this
map is one of the singularities. The new feature of our
improvement is that it is scaling invariant.

\

\noindent As an application of the previous statement we have that
if $\rho < 4 \pi (1+\a)$, then the above {\em center of mass}
cannot coincide with the singularity if the energy $I_\rho(u)$ is
sufficiently negative. The consequence of
the above fact (with a proper localization near each singularity)
is that low sublevels of the functional $I_\rho$ inherit some
topology from the surface $\Sig$ with a certain number of points
removed. Using these considerations, we are able to prove the
following result, where $G(\Sig)$ denotes the genus of $\Sig$.

\begin{thm}\label{t:mr} Suppose $\alpha_j \in (0,1]$ for all $j =
1, \dots, m$ and that $\rho \in (4 \pi, 8 \pi)$, with $\rho \neq 4
\pi (1 + \alpha_j)$ for all $j$. Denote by $J_\rho$ the subset of
$\{p_1, \dots, p_m\}$ for which $\rho < 4\pi(1+\a_i)$. Then if
$(G(\Sig), |J_\rho|) \neq (0,1)$ problem \eqref{eq:e0} has a
solution.
\end{thm}

\begin{rem}
For $G(\Sig) > 0$ this result has been proved in \cite{btcmp02},
Corollary 6. The case of the sphere is more delicate: in fact in
\cite{barlintar} it is shown (via a Pohozaev identity) that on the
standard sphere $(S^2,g_0)$ \eqref{eq:e0} has no solution for
$m=1$ and $\rho \in (4 \pi, 4 \pi (1+\a))$, $\a > 0$, which is
precisely the case $(G(\Sig), |J_\rho|) = (0,1)$. Therefore, our
condition $(G(\Sig), |J_\rho|) \neq (0,1)$ is somehow sharp.

\end{rem}

\noindent To prove Theorem \ref{t:mr} we use a min-max scheme
which is performed in detail in Section \ref{s:var}. First we show
that there exist test functions $\var_{\a,\l,x}$ (where $\a$ is a
suitable parameter) which satisfy $I_\rho(\var_{\a,\l,x}) \to -
\infty$ and $\mu_{\var_{\a,\l,x}} \rightharpoonup \d_x$ as $\l \to
- \infty$ whenever $x \in \Sig \setminus J_\rho$. Notice that, by
our assumptions, $\Sig \setminus J_\rho$ is a non contractible
set. We then consider continuous maps $\mathfrak{h}$ from the {\em
topological cone} over (a retraction of) $\Sig \setminus J_\rho$
into $H^1(\Sig)$, which coincide with the $\var_{\a,\l,x}$'s on
the boundary of the cone. The lower bound on the functional
described above shows that the supremum of $I_\rho$ on the image
of $\mathfrak{h}$ has a uniform control from below, provided $\l$
is sufficiently large. This allows us to show the admissibility of
the variational class consisting of the above maps $\mathfrak{h}$
(in this step the non contractibility of $\Sigma \setminus J_\rho$
is used), and to find Palais-Smale sequences for $I_\rho$ at some
bounded level.

At this point one can use a monotonicity result developed
initially by Struwe, which consists in varying the parameter
$\rho$ and to show that for a sequence $\rho_n \to \rho$ there
exist bounded Palais-Smale sequences and hence solutions to
\eqref{eq:e0}. The argument can then be completed using the a
priori estimates of \cite{btcmp02}, which imply compactness of
solutions for $\rho$ belonging to the ranges in Theorem
\ref{t:mr}.

\

\noindent The plan of the paper is the following. In Section
\ref{s:pr} we collect some preliminary results on the
Moser-Trudinger inequality (plus some more or less known
improvements), together with some compactness results and a
deformation lemma. In Section \ref{s:cov} we use a covering
argument to define a convenient center of mass for the measure
$\mu_u$ (see the above notation). In Section \ref{s:impr} we
obtain our new improved inequality, and lower bounds for $I_\rho$.
Section \ref{s:var} is devoted to the proof of Theorem \ref{t:mr};
some final comments are given at the end.

 \

\begin{center}

 {\bf Acknowledgements}

\end{center}

\noindent A.M. has been supported by the project FIRB-Ideas {\em
Analysis and Beyond}, and is grateful to the Mathematical Analysis
department of the University of Granada for the kind hospitality.
D.R has been supported by the Spanish Ministry of Science and
Innovation under Grant MTM2008-00988 and by J. Andaluc\'{\i}a (FQM
116). He thanks SISSA for the hospitality during a stay in 2009.
Both authors are grateful to D.Bartolucci and to the referee for some 
helpful comments.

\noindent

\section{Notation and preliminaries}\label{s:pr}

\noindent In this section we fix our notation and  recall some
useful known facts. We state in particular some variants and
improvements of the Moser-Trudinger inequality, together with
their consequences.

We  write $dist(x,y)$ to denote the distance between two points
$x, y \in \Sig$. Moreover, the symbol $B_p(r)$ stands for the open
metric ball of radius $r$ and center $p$, and $A_p(r,R)$ the
corresponding open annulus. $H^1(\Sig)$ is the Sobolev space of
the functions on $\Sig$ which are in $L^2(\Sig)$ together with
their first derivatives. The symbol $\| \cdot \|$ will denote the
norm of $H^1(\Sig)$. If $\Sig$ has boundary, $H^1_0(\Sig)$ will
denote the completion of $C_c^\infty(\Sig)$ with respect to the
Dirichlet norm. If $u \in H^1(\Sig)$, $\ov{u} = \frac{1}{|\Sig|}
\int_{\Sig} u dV_g = \intbar_\Sig u dV_g$ stands for the average
of $u$. For a real number $a$ we denote by $I_\rho^a$ the set $\{u
\in H^1(\Sig): \ I_\rho(u) \leq a\}$.

Large positive constants are always denoted by $C$, and the value
of $C$ is allowed to vary from formula to formula and also within
the same line. When we want to stress the dependence of the
constants on some parameter (or parameters), we add subscripts to
$C$, as $C_\d$, etc.. Also constants with this kind of subscripts
are allowed to vary.

\subsection{Improved Moser-Trudinger inequalities}

We start by recalling the well known Moser-Trudinger inequality in
a version that, when applied to the sphere, has also received the
name of Onofri inequality (see e.g. \cite{cygc}).

\begin{pro} Let $\Sig$ be a compact surface. Then
\noindent \begin{enumerate} \item[a)] If $\Sig$ has a boundary,
\begin{equation}\label{eq:mtdisk}
    \log \int_{\Sig} e^{2 u} dV_g \leq \frac{1}{2 \pi} \int_{\Sig}
  |\n_g u|^2 dV_g + 2 \fint_{\Sig} u + C \quad \hbox{ for every } u \in
  H^1(\Sig).
\end{equation}
\item[b)] If $\Sig$ has a boundary,
\begin{equation}\label{eq:mtdisk0}
    \log \int_{\Sig} e^{2 u} \ dV_g \leq \frac{1}{4 \pi} \int_{\Sig}
  |\n_g u|^2 dV_g + C \quad \hbox{ for every } u \in H^1_0(\Sig).
\end{equation}
\item[c)] If $\Sigma$ does not have a boundary, then
\begin{equation}\label{eq:mtdisk-1}
    \log \int_{\Sigma} e^{2 u} dV_g \leq \frac{1}{4 \pi} \int_{\Sig}
  |\n_g u|^2 dV_g + 2 \fint_{\Sig} u+ C \quad \hbox{ for every } u \in H^1(\Sig).
\end{equation}
\end{enumerate}
\end{pro}

\noindent The constant $\frac{1}{4\pi}$ in \eqref{eq:mtdisk-1} is
sharp, as on can see by using standard bubbles, peaked at some
point of $\Sigma$ (see \eqref{eq:bubble}). The constant in
\eqref{eq:mtdisk} is instead multiplied by two, since one can
center a bubble on the boundary of $\Sig$, dividing
approximatively by two both the conformal volume and the Dirichlet
energy. In \eqref{eq:mtdisk0} we impose $u\in H^1_0(\Sig)$, so
that this phenomenon is ruled out and the constant becomes again
$4 \pi$.

We begin by giving a {\em localized version} of the
Moser-Trudinger inequality,  following the ideas of \cite{dm}.

\begin{pro}\label{l:imprc}  Assume that $\Sig$ is a compact surface
(with or without boundary), and $\tilde{h}: \Sigma \to \R$
measurable, $0 \leq \tilde{h}(x) \leq C_0 \ a. e. \ x \in \Sig$.
Let $\Omega \subset \Sig$, $\delta>0$ such that
$dist(\Omega,\partial \Sig)>\delta$.

Then, for any $\e > 0$ there exists a constant $C = C(C_0, \e,
\delta)$ such that for all $u \in H^1(\Sigma)$,
$$
\log \int_{\Omega} \tilde{h}(x) e^{2 u } dV_g \leq \frac{1}{4 \pi -
\e} \int_{\Sig} |\n_g u|^2 dV_g + 2 \fint_{\Sig} u dV_g + C.
$$
\end{pro}

\begin{pf}
We can assume that $\fint_{\Sig} u=0$. Let us decompose
$$
u = u_1 + u_2,
$$
where $u_1 \in L^\infty(\Sigma)$ and $u_2 \in H^1(\Sigma)$ will be
fixed later. We have
\begin{equation}\label{eq:ddmm2}
    \log \int_{\Omega} \tilde{h}(x) e^{2 (u_1+u_2)} dV_g \leq 2
    \| u_1\|_{L^{\infty}(\Omega)} + \log \int_{\Omega}
\tilde{h}(x)e^{2 u_2} dV_g + C. \end{equation} We next consider a
smooth cutoff function $\chi$ with values into $[0,1]$ satisfying
$$
  \left\{
    \begin{array}{ll}
      \chi(x) = 1 & \hbox{ for } x \in \Omega,\\
      \chi(x) = 0 & \hbox{ if } dist(x, \Omega) > \delta/2,
    \end{array}
  \right.
$$
and then define
$$
  \tilde{u}(x) = \chi(x) u_2(x).
$$
Clearly, $\tilde{u} \in H_0^1(\Sig)$, so we can apply inequality
\eqref{eq:mtdisk0} to $\tilde{u}$, finding
$$
 \log \int_{\Omega} \tilde{h}(x) e^{2 u_2} dV_g \leq \log \int_{\Sigma}
 e^{2 \tilde{u}} dV_g + C \leq \frac{1}{4\pi}
  \int_{\Sigma} |\n (\chi(x) u_2(x))|^2 dV_g + C.
$$
Using the Leibnitz rule and the H\"older inequality we obtain:
$$
  \int_{\Sig} |\n (\chi(x) u_2(x))|^2 dV_g \leq (1+\e) \int_{\Sig}
|\n u_2|^2 dV_g + C_{\e} \int_{\Sig}
  |u_2|^2 dV_g.
$$
From \eqref{eq:ddmm2} and the last formulas we find
\begin{equation}\label{eq:last}
    \log \int_{\Omega} \tilde{h}(x) e^{2 u } dV_g \leq
    \frac{1+\e}{4\pi} \int_{\Sig}|\n u_2|^2 dV_g +
    C_{\e} \int_{\Sig} |u_2|^2 dV_g + 2 \|u_1\|_{L^\infty(\Omega)}
    +C.
\end{equation}
To control the latter terms we use truncations in Fourier modes.
Define $V_{\e}$ to be the direct sum of the eigenspaces of the
Laplacian on $\Sig$ (with Neumann boundary conditions) with
eigenvalues less or equal than $C_{\e}\e^{-1}$. Take now $u_1$ to
be the orthogonal projection of $u$ onto $V_{\e}$. In $V_\e$ the
$L^\infty$ norm is equivalent to the $L^2$ norm: by using
Poincar{\`e}'s inequality we get
$$ C_{\e} \int_{\Sig} |u_2|^2 dV_g \leq \e \int_{\Sig} |\n u_2|^2 dV_g, $$$$
  \|u_1\|_{L^\infty(\Omega)} \leq C'_{\e} \|u_1\|_{L^2(\Sig)} \leq C C'_\e \left(
  \int_{\Sig} |\n u_1|^2 dV_g \right)^{\frac 12} \leq
\e \int_{\Sig} |\n u_1|^2 dV_g +C''_{\e}.
$$
Hence, from \eqref{eq:last} and the above inequalities we derive
the conclusion by renaming $\e$ properly.
\end{pf}

\

\noindent The next result, for $\tilde{h}=1$ has been proved for
the first time in \cite{cl}. Assuming $\tilde{h}$ only bounded
does not require any changes in the arguments of the proof.
Roughly speaking, it states that if the function $e^{2u}$ is {\em
spread} into two regions of $\Sig$, then the constant in the
Moser-Trudinger inequality can be basically divided by $2$. Let us
point out that the arguments for Proposition \ref{l:imprc} could
also be used to prove the following result:

\begin{pro}\label{l:imprc0} Let $\Sig$ be a compact surface, $\tilde{h}:
\Sigma \to \R$ with $0 \leq \tilde{h}(x) \leq C_0$. Let $\Omega_1,
 \Omega_2$ be subsets of $\Sig$ with $dist(\Omega_1,\Omega_2)
\geq \delta_0$ for some $\delta_0>0$, and fix $\g_0 \in \left( 0,
\frac{1}{2} \right)$. Then, for any $\e > 0$ there exists a
constant $C = C(C_0, \e, \delta_0, \g_0)$ such that
$$
  \log \int_\Sig \tilde{h}(x)e^{2u} dV_g \leq C + \frac{1}{8
  \pi - \e} \int_{\Sig} |\n_g u|^2 dV_g + 2 \fint_{\Sigma} u.
$$
for all functions $u \in H^1(\Sig)$ satisfying
\begin{equation}\label{eq:ddmmi}
    \frac{\dis \int_{\Omega_i} \tilde{h}(x) e^{2u} dV_g}{\dis
\int_\Sig \tilde{h}(x) e^{2u} dV_g}
    \geq \g_0, \qquad \quad \quad  i =1,2.
\end{equation}

\end{pro}

%
%

\noindent A useful corollary of this result is the following one:
for the proof see \cite{cl} or \cite{djlw}.

\begin{cor}\label{c:conc} Suppose $\rho < 8 \pi$. Then, given any
$\e, r > 0$ there exists $L = L(\e,r) > 0$ such that
$$
  I_\rho(u) \leq - L \qquad \Rightarrow \qquad \dis \frac{\dis
  \int_{B_{x}(r)} \tilde{h} e^{2u} dV_g \,
  dV_g}{\dis \int_{\Sig} \tilde{h} e^{2u}\, d V_g} > 1 - \e \quad
  \hbox{ for some } x \in \Sig.
$$
\end{cor}

\subsection{Compactness of solutions and deformation lemma}

Concerning \eqref{eq:ee}, we have the following result, proved via
blow-up analysis.

\begin{thm}\label{t:bt2} (\cite{btcmp02})
Let $\Sig$ be a compact surface, and let $u_i$ solve \eqref{eq:ee}
with $\tilde{h}$ as in \eqref{eq:tildeh}, $\rho = \rho_i$, $\rho_i
\to \ov{\rho}$, with $\alpha_j
> 0$ and $p_j \in \Sig$. Suppose that $\int_\Sig
\tilde{h} e^{2 u_i} dV_g \leq C_1$ for some fixed $C_1
> 0$. Then along a subsequence $u_{i_k}$ one of the following
alternative holds:
\begin{description}
  \item[(i)] $u_{i_k}$ is uniformly bounded from above on $\Sig$;
  \item[(ii)] $\max_\Sig \left( 2 u_{i_k} - \log \int_\Sig \tilde{h} e^{2 u_{i_k}}
 dV_g  \right) \to + \infty$ and there exists a finite blow-up set
   $S = \{ q_1, \dots, q_l\} \in \Sig$ such that

   $(a)$ for any $s \in \{1, \dots, l\}$ there exist $x^s_n \to
   q_s$ such that $u_{i_k}(x^s_n) \to + \infty$ and $u_{i_k} \to
   - \infty$ uniformly on the compact sets of $\Sig \setminus
   S$,

   $(b)$ $\rho_{i_k} \frac{\tilde{h} e^{2 u_{i_k}}}{\int_\Sig
   \tilde{h} e^{2 u_{i_k}} d V_g} \rightharpoonup \sum_{s=1}^l
   \b_s \delta_{q_s}$ in the sense of measures, with $\b_s = 4
   \pi$ for $q_s \neq \{p_1, \dots, p_m\}$, or $\b_s = 4 \pi (1
   + \alpha_j)$ if $q_s = p_j$ for some $j = \{1, \dots, m\}$.
   In particular one has that
   $$
  \ov{\rho} = 4 \pi n + 4 \pi \sum_{j \in J} (1 + \alpha_j),
   $$
for some $n \in \N \cup 0$ and $J \subseteq \{1, \dots, m\}$
(possibly empty) satisfying $n + |J| > 0$, where $|J|$ is the
cardinality of the set $J$.
\end{description}
\end{thm}

\noindent From the above result we obtain immediately the
following corollary.

\begin{cor}\label{c:comp}
Suppose $\rho \in (4\pi, 8 \pi)$, and that
$$
  \rho \neq 4 \pi (1 + \alpha_j) \quad \hbox{ for all } j = 1, \dots, m.
$$
Then the set of solutions of \eqref{eq:ee} with zero mean value is
uniformly bounded in $C^2(\Sig)$.
\end{cor}

\

\noindent Our argument to prove existence of solutions relies on
variational theory:  more precisely, we look for some change in
the topology of the sublevels of $I_\rho$ which, via some
deformation lemma, leads to the existence of critical points. The
deformation lemma is usually employed when the {\em Palais-Smale
condition holds}: this means that every sequence $(u_l)_l$ for
which $I_\rho(u_l)$ converges and for which $I'_\rho(u_l)$ tends
to zero would admit a converging subsequence. This condition
allows to deform a sublevel into another, in case there are no
critical points in between, following the (negative) gradient flow
of the functional. The main role of the P-S condition is that the
flow lines stay compact as long as their energy is bounded.

Unfortunately it is still unknown whether the P-S condition holds
for $I_\rho$, and one has to bypass the argument via some other
kind of compactness result. We present next the following lemma,
obtained by M.Lucia (\cite{lucia}) through a variation of an
argument in \cite{str} (see also \cite{djlw}).

\begin{lem}\label{l:lucia} (\cite{lucia})
Given $a, b \in \R$, $a < b$, the following alternative holds:
either $\exists (\rho_l, u_l) \subseteq \R \times X$ satisfying
$$
  I'_{\rho_l}(u_l) = 0 \hbox{ for every } l; \qquad \qquad
  a \leq I_\rho(u_l) \leq b; \qquad \qquad \rho_l \to \rho,
$$
or the set $I_\rho^a$ is a deformation retract of $I_\rho^b$.
\end{lem}

\noindent In fact, the result in \cite{lucia} was proved for the
case of positive functions $\tilde{h}$, but it does extend to our
case as well with substantially the same proof.

\section{A covering argument}\label{s:cov}

\noindent In this section we use a covering argument (via thick
annuli) to detect both a {\em concentration size} for the function
$\tilde{h} e^{2u}$ and its location, which are useful to obtain
(in the next section) the desired improved inequality.

\begin{pro} \label{covering} Assume $\tilde{h}: \Sigma \to \R$,
$0 \leq \tilde{h}(x) \leq C$, $ \rho \in (4\pi, 8 \pi)$ and take a
constant $C_1>2$. There exist $\tau>0$, $L_0>0$ and a continuous
map
$$ \beta: I_{\rho}^{-L_0} \to \Sigma, $$
satisfying the following property: for any $u \in I_{\rho}^{-L_0}$
there exists $\bar{\sigma}
>0$ and $\bar{y} \in \Sigma$ such that $d(\bar{y}, \beta(u)) < 2
C_1 \bar{\sigma}$ and:
$$ \int_{B_{\bar{y}}(\bar{\sigma})} \tilde{h} e^{2u} dV_g
= \int_{\Sigma \setminus B_{\bar{y}}(C_1\bar{\sigma})} \tilde{h}
e^{2u} dV_g \geq \tau \int_{\Sigma} \tilde{h} e^{2u} dV_g.$$
\end{pro}

\begin{pf} Let us define:
$$ \mathcal{A}= \{ f \in L^1(\Sigma), \ f(x)>0 \ a.e., \ \int_{\Sigma} f dV_g
=1\}, $$
$$ \sigma: \Sigma \times \mathcal{A} \to (0,+\infty),  $$
where $\sigma=\sigma(x,f)$ is chosen such that:
$$ \int_{B_x(\sigma)} f dV_g = \int_{\Sigma \setminus B_x(C_1 \sigma)} f dV_g. $$
It is easy to check that $\sigma(x,f)$ is uniquely determined and
continuous. Moreover, we obtain that $\sigma$ satisfies
\begin{equation} \label{dett} dist(x,y) \leq C_1 \max \{ \sigma(x,f), \sigma(y,f)\}
+\min \{ \sigma(x,f), \sigma(y,f)\}.
\end{equation}
Otherwise, $B_x(C_1 \sigma(x,f)) \cap B_y(\sigma(y,f)+\e) =
\emptyset $ for some $\e>0$. Let us now show that
$A_y(\sigma(y,f), \sigma(y,f)+\e)$ is a nonempty open set.
Clearly, $B_y(\sigma(y,f)+\e)$ does not exhaust the whole surface
$\Sigma$. Since $\Sigma$ is connected, there exists $z \in
\partial B_{y}(\sigma(y,f)+\e)$. Then $dist(z,y)=
\sigma(y,f)+\e$, which  implies that $B_z(\e) \cap
B_{y}(\sigma(y,f)+\e)$ is a nonempty open set included in
$A_{y}(\sigma(y,f), \sigma(y,f)+\e)$.

Then:
$$ \int_{B_x(\sigma(x,f))} f dV_g =
\int_{\Sig \setminus B_x(C_1 \sigma(x,f))} f dV_g \geq
\int_{B_y(\sigma(y,f)+\e)} f dV_g > \int_{B_y(\sigma(y,f))} f
dV_g. $$ By interchanging the roles of $x$ and $y$, we would also
obtain the reverse inequality. This contradiction proves
\eqref{dett}.

Let us define $T: \Sigma \times \mathcal{A} \to (0,+\infty)$ by
$$T(x,f)= \int_{B_x(\sigma(x,f))} f dV_g.$$ Clearly, $T$ is also
continuous.

\bigskip

{\bf Step 1:} There exists $\tau>0$ such that $\max_{x\in \Sigma}
T(x,f) > 2 \tau $ for any $f \in \mathcal{A}$.

\bigskip
Let us take $x_0 \in \Sigma$ such that $T(x_0,f)=max_{x\in \Sigma}
T(x,f)$, and fix some $x \in A_{x_0}(\sigma(x_0,f ), C_1
\sigma(x_0,f))$.

We claim that:
\begin{equation} \label{desig1} dist(x,x_0)+ C_1 \sigma(x,f) \geq C_1
\sigma (x_0,f),\end{equation}
\begin{equation} \label{desig2} dist(x,x_0)- C_1 \sigma(x,f) \leq \sigma (x_0,f).
\end{equation}

Let us prove \eqref{desig1}. By contradiction, assume
$dist(x,x_0)+ C_1 \sigma(x,f) < C_1 \sigma (x_0,f)- 2\e$ for some
$\e>0$; by the triangular inequality, $B_x(C_1\sigma(x,f)) \subset
B_{x_0}(C_1 \sigma(x_0,f)-2\e)$. By  definition of $\sigma$,
$B_{x_0}(C_1 \sigma(x_0,f)) \neq \Sigma$. So we can show, as
previously, that $A_{x_0}(C_1 \sigma(x_0,f)-2\e, C_1
\sigma(x_0,f))$ is not empty. Then:
\begin{eqnarray*}
  T(x,f) & = & \int_{B_x(\sigma(x,f))} f dV_g =
\int_{\Sig \setminus B_x(C_1 \sigma(x,f))} f dV_g \\
   & > & \int_{\Sig \setminus B_{x_0}(C_1
\sigma(x_0,f))} f dV_g = T(x_0,f),
\end{eqnarray*}
which contradicts the definition of $x_0$.

We now prove \eqref{desig2} in an analogous way. Indeed, if
$dist(x,x_0)- C_1 \sigma(x,f) > \sigma (x_0,f)+2\e$, we obtain
that $\left (\Sig \setminus B(x, C_1 \sigma(x,f)) \right ) \supset
B(x_0, \sigma(x_0,f)+2\e)$. As above, the open set
$A_{x_0}(\sigma(x_0,f), \sigma(x_0,f)+\e)$ is nonempty. Then, we
obtain again a contradiction:
$$ T(x,f)= \int_{B_x(\sigma(x,f))} f dV_g =
\int_{\Sig \setminus B_x(C_1 \sigma(x,f))} f dV_g >
\int_{B_{x_0}(\sigma(x_0,f))} f dV_g = T(x_0,f). $$ The claim is
proved.

\bigskip

We subtract \eqref{desig2} from \eqref{desig1}, and deduce that
$\sigma(x,f) \geq \frac{C_1-1}{2C_1} \sigma (x_0,f) \geq
\frac{1}{4} \sigma(x_0,f)$ for any $x \in A_{x_0}(\sigma(x_0,f ),
C\sigma(x_0,f))$.

\medskip We now point out that given $C_1>2$, there exists $k=k(C_1)$ such that,
for any $\sigma >0$ and any $y \in \Sigma$,
$$ A_{y}(\sigma, C_1 \sigma) \subset \dis \cup_{i=1}^k
B(x_i, \frac 1 4 \sigma),$$ for some $x_i \in A_{y}(\sigma, C_1
\sigma)$.

\medskip Therefore:
$$ \int_{A_{x_0}(\sigma(x_0,f), C_1 \sigma(x_0,f))} f dV_g \leq \sum_{i=1}^k
\int_{B_{x_i}(\sigma(x_i,f))} f dV_g = \sum_{i=1}^k T(x_i,f) \leq k \,
T(x_0,f).$$ On the other hand:
$$ \int_{B_{x_0}(\sigma(x_0,f))} f dV_g = \int_{\Sig \setminus B_{x_0}(C_1
\sigma(x_0,f))} dV_g = T(x_0,f).$$ Hence $1= \int_{\Sigma} f dV_g
\leq (k+2) T(x_0,f)$, which concludes the proof of Step 1.

\bigskip

{\bf Step 2:} Definition of $\bar{\sigma}$ and $\bar{y}$.

\bigskip Let us define:
$$ S(f)= \{ x \in \Sigma:\ T(x,f) \geq \tau\}.$$

By Step 1, $S(f)$ is a nonempty compact set for every $f \in
\mathcal{A}$. Let us define also:
$$ \bar{\sigma}(f)= \max_{x \in S(f)} \sigma(x,f).$$
Observe that, in general, $\bar{\sigma} $ could be discontinuous
in $f$. Finally, take $\bar{y} \in S(f)$ such that
$\sigma(\bar{y},f)= \bar{\sigma}$.
For $u \in I_{\rho}^{-L}$, take $f= \left ( \int_{\Sigma}
\tilde{h} e^{2u} dV_g \right )^{-1} \tilde{h}(x) e^{2u(x)} =
\tilde{h}(x) e^{2u(x)-\log \int_{\Sigma} \tilde{h} e^{2u} dV_g}
\in \mathcal{A}$.

\bigskip

{\bf Step 3:} For any $\e>0$ there exists $L_0>0$ large enough
such that $diam \, S(f) \leq (C_1+1) \bar{\sigma} < \e$ for $u\in
I_{\rho}^{-L_0}$.

\bigskip

\noindent The first inequality holds independently of $L_0$;
indeed, by \eqref{dett}, $dist(x,y) \leq (C_1+1) \bar{\sigma}$ for
any given any $x$, $y \in S(f)$.

On the other hand, recall that:
$$\int_{B_{\bar{y}}(
\bar{\sigma})} \tilde{h} e^{2u} dV_g \geq \tau \int_{\Sigma} \tilde{h}
e^{2u} dV_g, \mbox{ and}$$
$$\int_{\Sig \setminus B_{\bar{y}}(
C_1\bar{\sigma})}  \tilde{h} e^{2u} dV_g \geq \tau \int_{\Sigma}
\tilde{h} e^{2u} dV_g.$$ Corollary \ref{c:conc} implies that we
can choose $L_0>0$ so that $\bar{\sigma}<\frac{\e}{C_1+1}$ for any
$u\in I_{\rho}^{-L_0}$.

\bigskip

{\bf Step 4:} Definition of $\beta(u)$ and conclusion.

\bigskip

\noindent We can assume to have an embedding $\Sigma \subset
\R^3$. Let us define: $$ \label{eta2} \eta: I_{\rho}^{-L_0} \to
\R^3, \ \ \eta(u)= \frac{\dis \int_{\Sigma} [T(x,f) - \tau]^+ x
dV_g}{\dis \int_{\Sigma} [T(x,f) - \tau]^+ dV_g} \ \ \mbox{ where
} f= \tilde{h}(x)e^{2u(x)-\log \int_{\Sigma} \tilde{h} e^{2u}
dV_g}. $$ Observe that in the above terms the integrands are equal
to zero outside $S(f)$.

Let $U \supset \Sigma$, $U \subset \R^3$ an open tubular
neighborhood of $\Sigma$, and $P: U \to \Sigma$ an orthogonal
projection onto $\Sigma$. Thanks to Step 3, $S(f) \subset
\overline{B}_{\bar{y}}( (C_1 +1) \bar{\sigma}) \subset
\overline{B}_{\bar{y}}^{\R^3}( (C_1 +1) \bar{\sigma})$. Since
$\eta(u)$ is a barycenter of a function supported in $S(f)$, we
have:
\begin{equation} \label{eta} |\eta(u) - \bar{y}| \leq (C_1 +1)\bar{\sigma}.\end{equation}
By taking $L_0>0$ large enough, $\eta(u) \in U$ for any $u \in
I_{\rho}^{-L_0}$. Therefore, we can define:
$$ \beta: I_{\rho}^{-L_0} \to \Sigma,\ \beta(u)= P \circ \eta(u).$$ To conclude the proof we just need to
show that $dist(\beta(u), \bar{y}) < 2C_1 \bar{\sigma}$. We denote
by $T_{\bar{y}} \Sigma$ the tangent space to $\Sigma$ at
$\bar{y}$. For any $x \in S(f) \subset B_{\bar{y}}((C_1 +1)
\bar{\sigma})$, one has:
\begin{equation} \label{hola1} \min \{ |\bar{y}+ y - x|: \ y \in T_{\bar{y}} \Sigma \} \leq C
\bar{\sigma}^2, \end{equation} where $C$ depends only on the $C^2$
regularity of $\Sigma$. Since $\eta(u)$ is a barycenter of a
function supported in $S(f)$, again, we have that
$$ \eta(u) \in \overline{B}_{\bar{y}}^{\R^3}((C_1 +1)
\bar{\sigma}), \ \min \{ |\bar{y}+ y - \eta(u)|: \ y \in
T_{\bar{y}} \Sigma \} \leq C \bar{\sigma}^2.$$ By taking a larger
$L_0$, if necessary, we can assume $2C \bar{\sigma}^2 \leq
\bar{\sigma}$ (recall, again, Step 3). So,
$$ |\beta(u) - \eta (u)| = \min_{x \in \Sigma} |\eta(u) - x| \leq 2C \bar{\sigma}^2
 \leq \bar{\sigma}. $$
This inequality, together with \eqref{eta}, implies that
\begin{equation} \label{beta} |\beta(u) - \bar{y}| \leq
(C_1 +2)\bar{\sigma}.\end{equation} Finally, take $\nu =
\frac{2C_1}{C_1+2}>1$; by Step 3 we can take $L_0$ larger, if
necessary, so that $\bar{\sigma}$ verifies that for any $x$, $y
\in \Sigma$, if $|x-y| \leq (C_1 +2)\bar{\sigma}$, then $dist(x,y)
\leq \nu |x-y|$. This, together with \eqref{beta}, finishes the
proof.
\end{pf}

\begin{rem}\label{r:beta} We claim that, with the above construction,
if $f_n=\frac{\tilde{h} e^{2 u_n}}{\int_\Sig \tilde{h} e^{2 u_n} dV_g} \rightharpoonup
\d_x$ for some $x \in \Sig$  then one also has $\beta(u_n) \to x$.

\medskip To check this, take $\bar{\sigma}_n = \bar{\sigma}(f_n)$, and $x_n \in
S(f_n)$. Passing to a subsequence, we have $\bar{\sigma}_n \to
\sigma_0$ and $x_n \to x_0$.

We first prove that $\bar{\sigma}_0 = 0$. If not, by construction,
we would have:

$$\int_{B(x_0, \sigma_0/2)} f_n dV_g > \tau,\ \ \int_{\Sig \setminus B(x_0, (C_1-1)\sigma_0)}
f_n dV_g > \tau,$$ which is a contradiction.

We now prove that $x_0=x$. If not, take $0<\delta <
dist(x_0,x)/2$. Since $\bar{\sigma}_n \to 0$,

$$\int_{B(x_0, \delta)} f_n dV_g > \tau,$$
which is again a contradiction.

Recall now that, by \eqref{beta}, $dist(\beta(u_n), S(f_n)) \leq C
\bar{\sigma}_n$; this concludes the proof of the claim.

\end{rem}

\begin{rem} \label{radial} If $\Sigma = B_0(1) \in \R^2$ and $u$
is a radial function, then also $f$ is a radial function.
Moreover, $\sigma(x,f)$ and $T(x,f)$ are radial functions, and the
set $S(f)$ is radially symmetric. Therefore, $\eta(u)=0$, and here
the projection $\beta$ coincides with $\eta$. So, any radial
function has zero barycenter.

This argument applies also to less restrictive symmetry
assumptions: for instance, if $u$ is even with respect to both the
$x$ and the $y$ axes.

Observe that in this case $\beta$ is defined in $H^1(B_0(1))$, and
not only on sublevels.

\end{rem}

\section{Improved inequalities}\label{s:impr}

\noindent The main goal of this section is to prove the following
result, giving a lower bound on the functional $I_\rho$ under
suitable conditions on its argument. Let $\beta$ be the map
constructed in Proposition \ref{covering}.

%

\begin{pro}\label{p:lb}
Assume $\rho \in (4 \pi, 4 \pi (1+\alpha_i))$, $C_1>1$
sufficiently large, $L_0>0$, $\tau>0$ such that Proposition
\ref{covering} applies. Then, there exists $L>L_0$ such that
$I_{\rho}(u)>-L$ for any $u \in I_{\rho}^{-L_0}$ satisfying that
$\beta(u)=p_i$.

\end{pro}

\begin{rem} \label{radial2} In particular, if $\Sigma=B_0(1)
\subset \R^2$, $p=0$ and $u$ is radially symmetric, we obtain the
boundedness below of $I_{\rho}$ for $\rho < 4 \pi (1+\alpha)$ (see
Remark \ref{radial}). This has already been observed in
\cite{det}.

The same thing is true under less restrictive symmetry
assumptions: for instance, for functions $u$ that are even with
respect to both the $x$ and the $y$ axes (see, again, Remark
\ref{radial}).
\end{rem}

\begin{pf}
Take $\delta>0$ fixed, $u \in I_{\rho}^{-L_0}$ such that
$\beta(u)=p_i$, and let $\bar{y} \in \Sigma$, $\bar{\sigma}>0$ be
as in Proposition \ref{covering}. For simplicity, let us assume
that $\int_{\Sigma} u dV_g =0$.

Take $\e>0$ (to be fixed later), and choose $s\in (2\bar{\sigma},
\frac{C_1}{2}
\bar{\sigma})$ such that: %
\begin{equation} \label{gamma}
\int_{A_{\bar{y}}(s/2,2s)} |\n u|^2 \, d V_g < \frac{1}{\gamma}
\int_{B_{\bar{y}}(\delta)} |\n u|^2\, d V_g, \ \mbox{ where }
\gamma \in \N,\ \ \e^{-1} < \gamma \leq \dis \frac{\log_2
C_1}{2}.\end{equation} Let us define:
$$ \cD_1= \int_{B_{\bar{y}}(s)} |\n u|^2\, d V_g, \quad
\cD_2= \int_{\Sig \setminus B_{\bar{y}}(s)} |\n u|^2\, d V_g,
$$$$\cD= \cD_1 + \cD_2, \quad \cJ= \log \int_{\Sigma} \tilde{h}(x) e^{2u} \, d V_g.
$$
Here the proof proceeds in three steps:

\medskip

\noindent {\bf Step 1:} We apply Proposition \ref{l:imprc} to
$u$ or, more precisely, to a convenient dilation of $u$ given by:
$$ v(x)= u(s x+ \bar{y}).$$
We have:
$$ \int_{B(\bar{y},s)} |\n u |^2 \, dV_g = \int_{B(0,1)} |\n v |^2 \, dV_g,  $$
$$ \fint_{B(\bar{y},s)} u  dV_g =  \fint_{B(0,1)} v dV_g,  $$
$$ \int_{B(\bar{y},s/2)} \tilde{h}(x) e^{2u} \, dV_g \leq C
\int_{B(\bar{y},s/2)} |x-p_i|^{2\alpha_i} e^{2u} \, dV_g \leq C
s^{2\alpha_i }\int_{B(\bar{y},s/2)} e^{2u}\, dV_g=$$$$ C
s^{2+2\alpha_i} \int_{B(0,1/2)} e^{2v} dV_g.$$ In the above
computations we have used that $|\bar{y} -p_i| \leq C s$. By
applying Proposition \ref{l:imprc} to $v$ and taking into account
the inequality $\int_{B_{\ov{y}}(s/2)} \tilde{h} e^{2u} \geq \tau
\int_{\Sig} \tilde{h} e^{2u}$ (recall Proposition \ref{covering}),
we obtain:

\begin{equation} \label{fund} \cJ - \frac{1}{4\pi-\e} \cD_1  + 2(1+\alpha_i) \log(1/s) \leq 2 \fint_{B(\bar{y},s)} u + C.
\end{equation}
This inequality is one of the key ingredients of the proof.

\medskip

\noindent {\bf Step 2:} We estimate $ \fint_{\partial
B_{\bar{y}}(s)} u$. To begin, consider a fixed value $s$ and the
function $\tilde{u}=u - \fint_{B_{\bar{y}}(s)} u$. By the trace
embedding, $\tilde{u} \in L^1(\partial B_{\bar{y}}(s))$, and
thanks to the Poincar{\`e}-Wirtinger inequality we get
$$ \Big| \fint_{\partial B_{\bar{y}}(s)}\tilde{u}\, dx \Big| \leq
C \| \tilde{u}\|_{H^1} \leq C \Big(\int_{B(\bar{y},s)} |\n u|^2\,
dV_g \Big)^{1/2}.$$ Therefore,
\begin{equation}
\label{xxx} \Big| \fint_{\partial B_{\bar{y}}(s)} u dV_g  -
\fint_{B_{\bar{y}}(s)} u dV_g \Big| \leq C
\Big(\int_{B(\bar{y},s)} |\n u|^2\, dV_g \Big)^{1/2} \leq \e
\cD_1+C'.
\end{equation} Observe now that the above inequality is invariant
under dilation. So, the constant $C$ is independent of $s$, and
hence $C'$ depends only on $\e$.

\medskip

\noindent {\bf Step 3: } By taking into account that $\tilde{h}(x)
\sim d(x,p_i)^{2\alpha_i}$ near $p_i$, and $|x-p_i| \leq C |x -
\ov{y}|$, we get the following estimate:
\begin{equation}\label{eq:stupid}
    \int_{\Sig \setminus B_{\bar{y}}(s)} \tilde{h}(x) e^{2u} dV_g =
    \int_{\Sig \setminus B_{\bar{y}}(s)} \frac{\tilde{h}(x)}{|x-\ov{y}|^{2 \a_i}} |x-\ov{y}|^{2 \a_i} e^{2u}
    dV_g \leq \frac{C}{s^{2\a_i}} \int_{\Sig \setminus B_{\bar{y}}(s)} e^{2v} dV_g,
\end{equation}
with $v(x) = \hat{u}(x) + 2 \a_i w(x)$,
$$w(x)= \left \{ \begin{array}{ll} \log s & x \in B_{\bar{y}}(s),
\\  \log |x| & x \in A_{\bar{y}}(s, \delta) \\ \log \delta & x \in \Sig
\setminus B_{\bar{y}}(\delta)
\end{array} \right. \, , \qquad
\left \{ \begin{array}{cc} -\Delta_g \hat{u}=0 \ \ x \in B(\bar{y},s), \\
\hat{u}(x)= u(x) \ \ x \notin B(\bar{y},s). \end{array} \right.
$$
Our intention is to apply the Moser-Trudinger to $v$. Observe
that:
$$ \fint_{\Sigma} v \leq C + \fint_{\Sigma} \hat{u}. $$
Since $\fint_{\Sig} u =0$ and $\hat{u} -u$ has compact support in
$B_{\bar{y}}(s)$,
$$ \left |\fint_{\Sigma} \hat{u}\right | = \left | \fint_{\Sigma} (\hat{u} -u) \right | \leq
C \left (\int_{B_{\bar{y}}(s)} |\n \hat{u} -\n u|^2 dV_g \right )^{1/2}
\leq \e \cD +C_{\e}. $$ We now estimate the Dirichlet energy:
$$ \int_{B_{\bar{y}}(s)} |\n v|^2 dV_g =\int_{B_{\bar{y}}(s)} |\nabla \hat{u}|^2 \, d V_g \leq C_0
\int_{A_{\bar{y}}(s/2,2s)} |\n u|^2 \, d V_g < C_0 \e \cD,$$ where
$C_0$ is a universal constant. Moreover,
\begin{eqnarray*}
  \int_{\Sig \setminus B_{\bar{y}}(s)} |\n v|^2 dV_g & = &
  \int_{\Sig \setminus B_{\bar{y}}(s)} |\n u|^2
dV_g + 4 \a_i^2 \int_{\Sig \setminus B_{\bar{y}}(s)} \frac{1}{|x-\ov{y}|^2} dV_g \\
  & + & 4 \a_i \int_{\Sig \setminus B_{\bar{y}}(s)} \n u \cdot \n
(\log|x-\ov{y}|) dV_g.
\end{eqnarray*}
We integrate by parts to obtain:
$$\int_{\Sig \setminus B_{\bar{y}}(s)} |\n v|^2 dV_g \leq \cD_2 + 8 \pi \a_i^2 \log \frac 1 s -
8 \pi \a_i \fint_{\partial B_{\bar{y}}(s)} u dV_g + C.$$ Next, we
apply the Moser-Trudinger inequality to $v$:
$$
\cJ \leq 2 \a_i \log \frac 1 s + \frac{1}{4\pi} \cD_2 + 2 \a_i^2
\log \frac 1 s - 2 \a_i \fint_{\partial B_{\bar{y}}(s)} u dV_g +
C_0 \e \cD+ C'.
$$
By using \eqref{fund} and \eqref{xxx}, we obtain:
$$
  (1 + \a_i) \cJ \leq \frac{\a_i}{4\pi-\e} \cD_1 + \frac{1}{4\pi} \cD_2 +C_0 \e \cD + \e \cD_1+ C''.
$$
In order to finish the proof it suffices to take $\e$ small enough
depending only on $\rho$ and $C_0$, (we recall that $C_0$ is a
universal constant). In such case, we need to choose $C_1$ large
enough so that \eqref{gamma} can be satisfied.
\end{pf}

\begin{rem}
The inequality used in \eqref{eq:stupid} is sharp when evaluated
on solutions of the  singular equation $- \D u = |x|^{2\a} e^{2u}$
in $\R^2$ (and properly glued on $\Sig$), which have been
classified in \cite{pratar} (see also \cite{chakie}) as
$$
  \var_{\a,\l}(y) = \log \frac{\l^{\a+1}}{\left( 1 +
  (\l |x|)^{2(1+\a)} \right)}.
$$
If $u$ has this form, and if $v = u + 2 \a \log |x|$, then $v$ has
the asymptotic profile of a standard bubble, which makes the
coefficients in the standard Moser-Trudinger inequality optimal.
\end{rem}

\

\noindent Let $J_\rho$ be as in Theorem \ref{t:mr}, and given a
small positive number $\th$ we define the compact set
\begin{equation}\label{eq:tThrho}
   \Theta_\rho = \Sig \setminus \cup_{p_i\in J_{\rho}} B_{p_i}(\th).
\end{equation}
From an extension of the above arguments one can prove the
following result.

\begin{pro}\label{p:proj} For $L > 0$ sufficiently large there
exists a continuous projection $\Psi$ from $I_\rho^{- L}$ into
$\Theta_\rho$ with the property that if $\frac{\tilde{h} e^{2
u_n}}{\int_\Sig \tilde{h} e^{2 u_n} dV_g} \rightharpoonup \d_x$
for some $x \in \Theta_\rho$  then $\Psi(u_n) \to x$.
\end{pro}

\begin{pf} It is sufficient to modify properly the function $\beta$ constructed
in Proposition \ref{covering}. Let us notice that since $\rho < 4
\pi (1 + \a_i)$ for $p_i \in J_\rho$, if $\b(u) = p_i$, $p_i \in
J_\rho$, then by Proposition \ref{p:lb} (and the subsequent
observation) $I_\rho$ is uniformly bonded from below. It follows
that if $L$ is sufficiently large and if $u \in I_\rho^{- L}$,
then $\beta(u) \in \Sig \setminus J_\rho$.

Then, if $\beta(u) \not\in \Theta_\rho$, it will belong to some
set of the form $B_{p_i}(\theta) \setminus \{ p_i \}$, for some
$p_i \in J_\rho$. At this point it is sufficient to move
$\beta(u)$ along the geodesic segment emanating from $p_i$ in the
direction of $\beta(u)$ until we hit the boundary of
$\Theta_\rho$. This procedure is well defined if $\theta$ is
chosen sufficiently small.

The last statement of the proposition follows from Remark
\ref{r:beta}.
\end{pf}

\section{Proof of Theorem \ref{t:mr} and final remarks}\label{s:var}

\noindent Let $\rho$, $J_\rho \subseteq \{p_1, \dots, p_m\}$ be as
in Theorem \ref{t:mr}, and let $\Theta_\rho$ be defined as in
\eqref{eq:tThrho}.

By our assumptions, the set $\Theta_\rho$ has either the topology
of $S^2$, or the topology of $S^2$ with more than one point
removed, or the topology of any other surface with $k$ points
removed, $k \geq 0$. In any of these situations, the set
$\Theta_\rho$ is non contractible.

We show next that it is possible to map an image of $\Theta_\rho$
homeomorphically into arbitrarily low sublevels of $I_\rho$. Let
$\tilde{\a} = \max_{i \in \{1, \dots, m\} \setminus J_\rho}
\alpha_i$. For $\l
> 0$ and $\a \in \left(\tilde{\a}, \frac{\rho}{4 \pi} - 1 \right)$,
we consider the following function
\begin{equation}\label{eq:varlx}
    \var_{\a,\l,x}(y) = \log \left( \frac{\l^{1+\a}}{1 +
    (\l dist(y,x))^{2(1+\a)}} \right).
\end{equation}

\begin{lem}\label{l:lowen}
Let $\var_{\a,\l,x}$ be as in \eqref{eq:varlx}. Then one has
$$
  I_\rho(\var_{\a,\l,x}) \to - \infty \qquad \hbox{ as } \l \to + \infty
$$
uniformly on $\Theta_\rho$.
\end{lem}

\begin{pf}
It is standard (see e.g. Section 4 in \cite{maldan}) to check that
\begin{equation}\label{eq:estbau}
  \int_\Sig |\n \var_{\a,\l,x}|^2 d V_g \simeq 8 \pi (1+\a)^2 \log \l;
  \qquad \quad \intbar_\Sig \var_{\a,\l,x} d V_g \simeq - (1 + \a) \log \l
\end{equation}
as $\l \to + \infty$. We want to estimate next the integral of the
exponential term. If $x$ belongs to a compact set of $\Sigma
\setminus \{p_1, \dots, p_m\}$ then for some $\d > 0$ we have that
\begin{equation}\label{eq:miao}
    \int_\Sig \tilde{h} e^{2 \var_{\a,\l,x}} dV_g \geq \frac{1}{C_\d}
    \int_{B_x(\d)} e^{2 \var_{\a,\l,x}} dV_g.
\end{equation}
Using normal geodesic coordinates $y$ centered at $x$ and the
estimate
\begin{equation}\label{eq:vold}
    dV_g = (1 + o_\d(1)) dy; \qquad d(x, y) = (1+o_\d(1)) |y|, \qquad
    y \in B_x(\d)
\end{equation}
(where we identified $y$ with its coordinate in the above system)
we find that
$$
   \int_{B_x(\d)} e^{2 \var_{\a,\l,x}} dV_g = (1+o_\d(1)) \int_{B_0(\d)}
   \frac{\l^{2(1+\a)}}{\left(1 + (\l |y|)^{2(1+\a)} \right)^2} dy
   \geq \frac{1}{C} \l^{2\a}.
$$
The last formula, with \eqref{eq:estbau}, \eqref{eq:miao} and the
fact that $\rho > 4\pi (1+\a)$ imply the conclusion for $x$ in
compact sets of $\Sigma \setminus \{p_1, \dots, p_m\}$.

Next, it is sufficient to consider the case in which $x$ is close
to one of the $p_i$'s with $i \not\in J_\rho$. Localizing the
integral as in \eqref{eq:miao} and \eqref{eq:vold}, it is
sufficient to estimate from below the following quantity
$$
  \int_{B_{p_i}(\d)} |y - p_i|^{2 \alpha_i} \frac{\l^{2(1+\a)}}{\left(
  1 + (\l |y-x|)^{2(1+\a)} \right)^2} dy,
$$
uniformly for $p_i \in \{p_1, \dots, p_m\} \setminus J_\rho$, and
for $|x - p_i| \leq \d^2$. By a change of variables we are left
with
$$
  \int_{B_{0}(\d)} \frac{\l^{2(1+\a)} |y|^{2 \alpha_i}}{\left( 1 +
  (\l |y-x|)^{2(1+\a)} \right)^2} dy; \qquad \quad |x| \leq \d^2.
$$
We next divide the domain into the three sets
$$
  \mathcal{B}_1 = \{ |y| \leq \sqrt{\d} |x| \}; \qquad
  \mathcal{B}_2 = \{ \sqrt{\d} |x| < |y| \leq \d^{-\frac 12} |x| \};
$$
$$
  \mathcal{B}_3 = \{ \d^{-\frac 12} |x| < |y| \leq \d \}.
$$
In $\mathcal{B}_1$ we have that $|y - x| = (1 + o_\d(1))|x|$,
which implies
\begin{eqnarray*}
  \int_{\mathcal{B}_1} \frac{\l^{2(1+\a)} |y|^{2 \alpha_i}}{\left(
  1 + (\l |y-x|)^{2(1+\a)} \right)^2} dy  & \geq & C_\d^{-1} \frac{\l^{2(1+\a)}}{\left(
  1 + (\l |x|)^{2(1+\a)} \right)^2} \int_{B_{0}(\sqrt{\d}|x|)} |y|^{2 \alpha_i} dy \\
   & \geq & C_\d^{-1} \frac{\l^{2(1+\a)} |x|^{2 \alpha_i + 2}}{\left(
  1 + (\l |x|)^{2(1+\a)} \right)^2}.
\end{eqnarray*}
In $\mathcal{B}_2$ $|y|$ is bounded above and below by constants
(depending on $\d$) multiplying $|x|$, and hence with a change of
variables we get
\begin{eqnarray*}
  \int_{\mathcal{B}_2} \frac{\l^{2(1+\a)} |y|^{2 \alpha_i}}{\left(
  1 + (\l |y-x|)^{2(1+\a)} \right)^2} dy & \geq & C_\d^{-1}
 \int_{\l(\mathcal{B}_2 - x)} \frac{|x|^{2 \alpha_i}
  \l^{2\a} dz}{(1 + |z|^{2(1+\a)})^2}  \\
  & \geq & C_\d^{-1} |x|^{2 \alpha_i} \l^{2\a} \int_0^{\frac{\l |x|}{C_\d}}
  \frac{s ds}{(1 + s^{1+\a})^2}.
\end{eqnarray*}
Integrating we obtain
$$
  \int_{\mathcal{B}_2} \frac{\l^{2(1+\a)} |y|^{2 \alpha_i}}{\left(
  1 + (\l |y-x|)^{2(1+\a)} \right)^2} dy \geq C_\d^{-1} |x|^{2 \alpha_i}
  \l^{2\a} \frac{\l^2 |x|^2}{1 + \l^2 |x|^2}.
$$
Finally in $\mathcal{B}_3$ we have that $|y-x|^2 = (1 +
o_\d(1))|y|^2$, and therefore from a change of variables we get
\begin{eqnarray*}
  \int_{\mathcal{B}_3} \frac{\l^{2(1+\a)} |y|^{2 \alpha_i}}{\left(
  1 + (\l |y-x|)^{2(1+\a)} \right)^2} dy & \geq & C_\d^{-1}
  \int_{\mathcal{B}_3} \frac{\l^{2(1+\a)} |y|^{2 \alpha_i} dy}{(1 +
  (\l |y|)^{2(1+\a)}))^2} \\
  & \geq & C_\d^{-1} \l^{2 (\a-\alpha_i)} \int_{C_\d \l |x|}^{\frac{\l}{C_\d}}
 \frac{s^{2 \alpha_i + 1} ds}{(1 + s^{2(1+\a)}))^2}.
\end{eqnarray*}
Evaluating the last integral we have that
$$
  \int_{\mathcal{B}_3} \frac{\l^{2(1+\a)} |y|^{2 \alpha_i}}{\left(
  1 + (\l |y-x|)^{2(1+\a)} \right)^2} dy \geq C_\d^{-1} \l^{2
  (\a-\alpha_i)} \frac{1}{1 + (\l |x|)^{2 + 4\a - 2 \alpha_i}}.
$$
In conclusion we obtain the estimate
$$
  \int_{B_{p_i}(\d)}  \frac{|y - p_i|^{2 \alpha_i} \l^2}{\left(
  1 + \l^2 |y-x|^2 \right)^2} dy \geq C_\d^{-1} \left( \frac{\l^{2(1+\a)}
  |x|^{2(1+\alpha_i)}}{1 + \l^2 |x|^2} + \frac{\l^{2
  (\a-\alpha_i)}}{1 + (\l |x|)^{2 + 4\a - 2 \alpha_i}} \right).
$$
Using elementary arguments, one can easily check
that the last quantity is always greater or equal to $C_\d^{-1}
\l^{2\a-2\alpha_i}$.

Therefore, adding the three terms in $I_\rho$ we find
\begin{eqnarray*}
  I_\rho(\var_{\a,\l,x}) & \leq & 8 \pi (1+\a)^2 \log \l - 2
  \rho (1+\a) \log \l - 2 \rho (\a - \alpha_i) \log \a + l.o.t. \\
   & = & 2 \log \l \left( 4 \pi (1+\a)^2
  - (1+\a) \rho - \rho (\a - \alpha_i) \right) + l.o.t.
\end{eqnarray*}
since $\rho > 4 \pi (1+\alpha_i)$ and since $\a > \alpha_i$, we
get the conclusion.
\end{pf}

\

\noindent We also have the following result, regarding the
concentration properties of the functions $\tilde{h} e^{2
\var_{\a,\l,x}}$.

\begin{lem}\label{l:conc} Let $\var_{\a,\l,x}$ be defined as in Lemma
\ref{l:lowen}. Then, for any $x \in \Theta_\rho$,
$$
  \frac{\tilde{h} e^{2 \var_{\a,\l,x}}}{\int_\Sig \tilde{h} e^{2 \var_{\a,\l,x}}
  dV_g} \rightharpoonup \d_x  \qquad \hbox{ as } \l \to +
  \infty.
$$
\end{lem}

\begin{pf}
Given $\e > 0$ it is sufficient to show that
$$
  \frac{\int_{\Sig \setminus B_{x}(\e)} \tilde{h} e^{2 \var_{\a,\l,x}}
  dV_g}{\int_\Sig \tilde{h} e^{2 \var_{\a,\l,x}} dV_g} \to 0
  \qquad \quad \hbox{ as } \l \to + \infty.
$$
By the proof of Lemma \ref{l:lowen} we derived that $\int_\Sig
\tilde{h} e^{2 \var_{\a,\l,x}} dV_g \geq C^{-1} \l^{2\a - 2
\hat{\a}}$, and therefore it is sufficient to check that
\begin{equation}\label{eq:chch}
    \frac{\int_{\Sig \setminus B_{x}(\e)} \tilde{h} e^{2 \var_{\a,\l,x}}
  dV_g}{\l^{2\a - 2 \hat{\a}}} \to 0
  \qquad \quad \hbox{ as } \l \to + \infty.
\end{equation}
For doing this, we notice that
$$
  e^{2 \var_{\a,\l,x}} \leq C_\e \l^{-2(1+\a)} \qquad \quad \hbox{ as } \l \to + \infty,
$$
and therefore \eqref{eq:chch} follows immediately.
\end{pf}

\

%
%

\noindent We next define the min-max scheme which is needed to
prove Theorem \ref{t:mr}. We fix $L>0$ as in Proposition
\ref{p:lb}, and then $\l > 0$ so large that
$I_\rho(\var_{\a,\l,x}) < - 2 L$ for $x \in \Theta_\rho$. The
latter choice is possible in view of Lemma \ref{l:lowen}.

We then define the set
$$
  \L_\l = \left\{ \var_{\a,\l,x} \; : \; x \in \Theta_\rho \right\}.
$$
Next, we consider the topological cone
$$
  \hat{\Theta}_\rho = \left( \Theta_\rho \times [0,1] \right)/(\Theta_\rho
  \times \{1\}),
$$
where the equivalence relation identifies all the points in
$\Theta_\rho \times \{1\}$. Let us introduce next the family of
continuous maps
$$
  \mathcal{H}_{\l,\rho} = \left\{ \mathfrak{h} : \hat{\Theta}_\rho \to H^1(\Sig) \; : \;
  \mathfrak{h}(x) = \var_{\a,\l,x} \hbox{ for every } x \in \Theta_\rho \right\},
$$
and the number
$$
  \ov{\mathcal{H}}_{\l,\rho} = \inf_{\mathfrak{h} \in \mathcal{H}_\l} \sup_{z \in \hat{\Theta}_\rho}
  I_\rho(\mathfrak{h}(z)).
$$
We have then the following result.

\begin{pro}\label{p:crit} Under the assumptions of Theorem \ref{t:mr},
if $\l$ is sufficiently large the number
$\ov{\mathcal{H}}_{\l,\rho}$ is finite. Moreover
$\ov{\mathcal{H}}_{\l,\rho}$ is a critical value of $I_\rho$.
\end{pro}

\begin{pf} If $L$ is as in Proposition \ref{p:lb}, we show that
indeed $\ov{\mathcal{H}}_{\l,\rho} > - \frac 32 L$. In fact,
suppose by contradiction that there exists a map $\mathfrak{h}_0$
such that
\begin{equation}\label{eq:h0}
 \mathfrak{h}_0 \in \mathcal{H}_{\l,\rho} \qquad \hbox{ and } \qquad
 \sup_{z \in \hat{\Theta}_\rho} I_\rho(\mathfrak{h}_0(z)) \leq - \frac 32 L.
\end{equation}
Then Proposition \ref{p:proj} applies and gives a continuous map
$F_{\l,\rho} : \hat{\Theta}_\rho \to \Theta_\rho$ defined as the
composition
$$
  F_{\l,\rho} = \Psi \circ \mathfrak{h}_0.
$$
since $\mathfrak{h}_0 \in \mathcal{H}_{\l,\rho}$, and hence it
coincides with $\var_{\a,\l,\cdot}$ on $\pa \hat{\Theta}_\rho
\simeq \Theta_\rho$, by Lemma \ref{l:conc} and Remark \ref{r:beta}
we deduce that
\begin{equation}\label{eq:hom}
      F_{\l,\rho}|_{\Theta_\rho}  \hbox{ is homotopic to } Id|_{\Theta_\rho}.
\end{equation}
Here, the homotopy is given by the parameter $\l$ as $\l \to
+\infty$.

\medskip Let us write as pairs $(x,\omega)$ the elements of
$\hat{\Theta}_\rho$. If we let $\omega$ run between $1$ and $0$,
and we consider the maps $F_{\l,\rho}(\cdot,\omega) : \Theta_\rho
\to \Theta_\rho$, we obtain an homotopy between
$F_{\l,\rho}|_{\Theta_\rho}$ and a constant map. Since by our
assumptions $\Theta_\rho$ is not contractible, we obtain a
contradiction with \eqref{eq:hom}. This proves the first part of
the statement.

\

\noindent To check that $\ov{\mathcal{H}}_{\l,\rho}$ is a critical
level, we use a monotonicity method introduced by Struwe, and
which has been used extensively in the study of Liouville type
equations. We consider a sequence $\rho_n \to \rho$ and the
corresponding functionals $I_{\rho_n}$. All the above estimates
and results can be worked out for $I_{\rho_n}$ as well with minor
changes.

We then define the number $\tilde{\mathcal{H}}_{\l,\rho} :=
\frac{\ov{\mathcal{H}}_{\l,\rho}}{\rho}$, which corresponds to the
functional $\frac{I_\rho}{\rho}$. It is immediate to see that
$$
  \rho \mapsto \tilde{\mathcal{H}}_{\l,\rho} \qquad \hbox{ is monotone },
$$
and, reasoning as in \cite{djlw}, there exists a subsequence of
$(\rho_n)_n$ such that $I_{\rho_n}$ has a solution $u_n$ at level
$\ov{\mathcal{H}}_{\l,\rho_n}$. Then, applying Theorem \ref{t:bt2}
and passing to a further subsequence, we obtain that $u_n$
converges to a critical point $u$ of $I_\rho$ at level
$\ov{\mathcal{H}}_{\l,\rho}$.
\end{pf}

\begin{rem}\label{r:ab} (a) In \cite{cl3} the Leray-Schauder degree of
\eqref{eq:ee} is being computed, using refined blow-up analysis
and Lyapunov-Schmidt reductions, in the spirit of \cite{cl1},
\cite{cl2}. Taking into account the result in \cite{mal}, we
speculate that there should be a relation between the degree and
the Euler characteristic of the set $\Sigma \setminus J_{\rho}$.
Anyway, the min-max methods might give existence (or even
multiplicity, see \cite{dem}, \cite{dem2}) results even when the
total degree of the equation vanishes. On the other hand, we have
to avoid some {\em critical} values of $\rho$ which are instead
treatable via blow-up analysis.

(b) Further results are discussed in \cite{bdm}, where the case of
arbitrary positive $\alpha$'s is discussed for  surfaces with
positive genus. In this case one can exploit the non simply
connectedness of the surfaces and avoid to use of the improved
inequality we derived here. However the latter should be necessary
to treat the case with singularities of different signs (or to
characterize the homology of low sublevels of $I_\rho$, as in
\cite{mald}). The case of the sphere for larger values of $\rho$
will be treated in a forthcoming paper, combining our approach
with some techniques in \cite{dm} and a topological argument. In
\cite{ barmon} and \cite{carl} the case of negative $\alpha$'s is
studied: in this situation one
can combine Troyanov's result with Chen-Li's inequality.
\end{rem}

\end{document}